\documentclass[twocolumn]{autart}    % Enable this line and disable the 
\usepackage{graphicx}          % Include this line if your 
\usepackage{natbib}
\usepackage{amssymb}
\usepackage{amsmath}
\usepackage{graphicx}
\graphicspath{ {fig/} }
\usepackage{epstopdf}
\usepackage{xcolor}
\usepackage{float}
\floatstyle{ruled}
\newfloat{algorithm}{thp}{lop}
\floatname{algorithm}{Algorithm}
\newfloat{subroutine}{thp}{lop}
\floatname{subroutine}{Subroutine}
\usepackage{xcolor}
\usepackage{booktabs} % For better table formatting
\usepackage{diagbox} % put in preamble

\renewcommand{\theenumi}{\roman{enumi})}%

\newtheorem{proposition}{Proposition}
\newtheorem{theorem}{Theorem}

\begin{document}
\begin{frontmatter}

\title{Projection-based discrete-time consensus on the unit sphere\thanksref{footnoteinfo}}

\thanks[footnoteinfo]{This paper was not presented at any IFAC 
meeting. Corresponding author is Johan Thunberg.}

\author[Lund]{Johan Thunberg}\ead{johan.thunberg@eit.lth.se},     
\author[Lund]{Galina Sidorenko}\ead{galina.sidorenko@eit.lth.se}   

\address[Lund]{Department of Electrical and Information Technology, Faculty of Engineering, Lund University}  
\maketitle


\begin{abstract}                          % Abstract of not more 
We address discrete-time consensus on the Euclidean unit sphere. For this purpose we consider a distributed algorithm comprising the  iterative projection of a conical combination of neighboring states. Neighborhoods are represented by a strongly connected directed graph, and the conical combinations are represented by a (non-negative) weight matrix with a zero structure corresponding to the graph. A first result mirrors earlier results for gradient flows. Under the assumptions that each diagonal element of the weight matrix is more than $\sqrt{2}$ larger than the sum of the other elements in the corresponding row, the sphere dimension is greater or equal to $2$, and the graph, as well as the weight matrix, is symmetric, we show that the algorithm comprises gradient ascent, stable fixed points are consensus points, and the set of initial points for which the algorithm converges to a non-consensus fixed point has measure zero. The second result is that for the unit circle and a strongly connected graph or for any unit sphere with dimension greater than or equal to $1$ and the complete graph, only for a measure zero set of weight matrices there are fixed points for the algorithm which do not have consensus or antipodal configurations.
\end{abstract}

\end{frontmatter}

\section{Introduction}
The consensus or synchronization problem, where agents with limited communication must agree on a shared quantity, such as position or velocity, has garnered significant attention across diverse fields, from engineering ~\cite{schenato2011average,savazzi2020federated,olfati2007consensus} to physics \cite{vicsek1995novel} and biology~\cite{cao2013overview,ren2005survey,reynolds1987flocks}. Many real-world applications and models involve agents evolving on nonlinear manifolds, including spheres~\cite{sepulchre2010consensus,sarlette2009synchronization} and rotation groups~\cite{tron2012intrinsic,thunberg2018dynamic}. 

This paper analyzes a distributed consensus algorithm for the unit sphere, where in each iteration conical combinations of neighboring states are projected onto the unit sphere. The neighborhoods are represented by a strongly connected directed graph, and the conical combinations are determined by a non-negative weight matrix with the same zero structure as the adjacency matrix of the graph. As such, the algorithm is perhaps one of the most natural choices in this context, and for symmetric graphs and symmetric weight matrices the algorithm comprises projected gradient ascent (or descent, depending on the choice of sign).

In the case of symmetric graphs and weight matrices, the problem falls into a large class of synchronization problems, dating back to at least the Kuramoto model~\cite{kuramoto1975self}, which has recently been investigated from the optimization perspective for unit spheres~\cite{geshkovski2023mathematical,townsend2020dense} and other manifolds~\cite{mcrae2024benign}. The unit circle, perhaps surprisingly, presents a greater challenge than higher-dimensional unit spheres when studying consensus/synchronization in the context of symmetric graphs and weight matrices~\cite{johanl2017continious}. 

Synchronization on the unit circle, $\mathbb{S}^1$, was well-summarized in~\cite{dorfler2014synchronization}, where, among other techniques, potential landscape analysis was reviewed. In such approaches the continuous-time synchronization protocol is the gradient of a potential function, and the Jacobian thereof corresponds to the Hessian of the potential function. Certain undirected connected graphs such as acyclic graphs, complete graphs, and sufficiently dense graphs, are called $\mathbb{S}^1$-synchronizing if all critical points not in the consensus set (where all the states are equal) are local maxima or saddle points~\cite{canale2008almost,canale2010complexity,canale2010wheels,monzon2005global,sarlette2009geometry,taylor2012there}. 

The concept of $\mathbb{S}^1$-synchronizing graphs was generalized to $\mathbb{S}^{d-1}$-synchronizing graphs in~\cite{johanl2017continious} for $d \geq 2$ being the dimension of the ambient space. Interestingly, for $d \geq 3$, all symmetric and connected graphs are $\mathbb{S}^{d-1}$ synchronizing, whereas for $\mathbb{S}^1$ this is not true. 
More general results for graphs that synchronize on Stiefel manifolds were later presented in~\cite{markdahl2020high}. These continuous-time protocols for undirected graphs do not rely on global parameters or shaping-functions, which was, for example, the case in the discrete-time approach in~\cite{tron2012intrinsic} for synchronization on $SO(3)$. The first results of this paper are similar to a subset of those in~\cite{johanl2017continious} but for discrete time, where an additional requirement is included for the weight matrix, see contribution i) below. 

As discussed, in contrast to the higher-dimensional spheres, for the unit circle there might be fixed points that are neither unstable nor have trivial consensus or antipodal configurations. The second and main contribution of this paper (see contribution ii) below) is to flip the perspective in this context and, instead of investigating the properties of unfavorable fixed points, investigate the occurrence of such fixed points over the space of weight matrices. It turns out that the probability is $0$ of randomly choosing a weight matrix corresponding to a chosen strongly-connected graph topology, for which there are fixed points neither in a consensus configuration nor in an antipodal configuration.   

The main contributions of this paper are as follows.
\begin{enumerate}
    \item A condition for similar results as presented for gradient flows in~\cite{johanl2017continious} for the special case of symmetric weights which are constant. For unit sphere dimension greater or equal to $2$, or ambient space dimension $d$ greater or equal to $3$, if the weight matrix is such that each diagonal element is more than $\sqrt{2}$ larger than the sum of the other elements in the corresponding row, stable fixed points are consensus points, and the set of initial points for which the algorithm converges to a non-consensus fixed point has measure zero. \\
    \item For the unit circle and strongly connected graph, or for any unit sphere with dimension greater than or equal to $1$ and the complete graph, only for a measure zero set of weight matrices there are fixed points for the algorithm which do not satisfy that state vectors are equal up to sign (consensus or antipodal configurations).
\end{enumerate}

The paper is structured as follows: Section \ref{sec:preliminaries} provides preliminary  notations and definitions. Section \ref{sec:algo} introduces our proposed algorithm along with properties thereof. Section \ref{sec:properties_diff} establishes a result on invertibility of the differential needed for contribution i) in Section \ref{sec:symmetric}. Furthermore, Section 4 contains the result that points in the consensus set are not unstable. Section \ref{sec:general_weight_matricec} provides the results for contribution ii). Section \ref{sec:simulations} presents simulation results mainly focused on contribution ii). Finally, Section \ref{sec:conclusions} summarizes the main findings and conclusions.

\section{Preliminaries} \label{sec:preliminaries}
 A directed graph with nodes \( \mathcal{V}(n) = \{1, \dots, n\} \) and edges \( \mathcal{E} \subset \mathcal{V}(n) \times \mathcal{V}(n) \) is denoted by \( \mathcal{G}(n) = (\mathcal{V}(n), \mathcal{E}) \). Such a graph is strongly connected if there exists a directed path between any two nodes \(i\) and \(j\). Such a graph is symmetric if $(i,j) \in \mathcal{E}$ implies $(j,i) \in \mathcal{E}$.

Throughout, vectors are taken to be column vectors by convention. For a matrix \( A \), \( [A]_i \) denotes the \( i \)-th row, \( [A]^j \) denotes the \( j \)-th column, and \( [A]_{ij} \) refers to the element in the \( i \)-th row and \( j \)-th column. If $a_{ij}$ is the element of $A$ at row $i$ and column $j$, we write $A = [a_{ij}]$. A matrix \( A = [a_{ij}] \in \mathbb{R}^{n \times n} \) is a weight matrix for \( \mathcal{G}(n) \) if all elements in the matrix are non-negative, and \( a_{ij} > 0 \) if and only if $(i,j) \in \mathcal{E}$ for all $i,j$ such that $i \neq j$. For a symmetric graph, the weight matrix is symmetric if $a_{ij}=a_{ji}$ $\forall i,j$. We say that a weight matrix \( A \) is strictly diagonally dominant if \( [A]_{ii} > \sum_{j \neq i} [A]_{ij} \)  \(\forall  i \). 

A column vector of ones in $\mathbb{R}^{n}$ is denoted by $\textbf{1}_n$. The identity matrix of size $n$ is denoted by $I_n$. The symbol $\otimes$ represents {the} Kronecker product of two matrices, and $\text{vec}(\cdot)$ is the the vectorization operation, which stacks the columns of a matrix consecutively into a single column vector. Finally, $\text{tr}(A)$ denotes the trace of a matrix $A$, and $\text{diag}\left (d_1, \dots d_n \right )$ denotes the diagonal matrix of size $n \times n$ with diagonal entries $d_1, \dots d_n$ and zeros elsewhere.

For $d \geq 2$, we denote the $(d-1)$-dimensional Euclidean unit sphere in $\mathbb{R}^d$ by
\begin{align}
\mathbb{S}^{d-1} = \left\{ x \in \mathbb{R}^d : \|x\|_2 = 1 \right\}.
\end{align}
We further define the set
\begin{align}
\mathbb{S}(n,d)\!=\!\{X \in \mathbb{R}^{n \times d}: \|[X]_i\|_2 = 1 \; \forall  i \in \{1,\dots, n\}\},
\end{align}
which is equivalent to $(\mathbb{S}^{d-1})^n$ below and comprises matrices in $\mathbb{R}^{n \times d}$ whose rows all have unit norm, i.e., each row is an element of $\mathbb{S}^{d-1}$. 

Furthermore, let 
\begin{align}
(\mathbb{S}^{d-1})^n = \{x \in \mathbb{R}^{nd}: \exists X \in \mathbb{S}(n,d) \text{ s.t. } x = \text{vec}(X^T)\}.
\end{align}

\section{Consensus algorithm} \label{sec:algo}
Let $A = [a_{ij}]$ be a strictly diagonally dominant weight matrix for a strongly connected directed graph $\mathcal{G}(n) = (\mathcal{V}(n), \mathcal{E})$. Let $x_i(k) \in \mathbb{S}^{d-1}$, for $d \geq 2$, $i \in \{1, \dots, n\}$, $k = 0,1,2,\ldots$, where 
\begin{align}\label{eq:main_algo}
x_i(k+1) = \frac{([A]_i \otimes I_d)x(k)}{\|([A]_i \otimes I_d)x(k)\|_2} = \frac{\sum_{j = 1}^n a_{ij}x_j(k)}{\|\sum_{j = 1}^n a_{ij}x_j(k)\|_2}.
\end{align}
%for $i = 1,2, \ldots, n$ and $k$ being non-negative integer. 
It is easy to verify that $A$ being strictly diagonally dominant ensures that $x_i(k)$ is well defined for all $i$ and $k$ in the sense that there are no divisions by zero. We will assume throughout this paper that $A$ is strictly diagonally dominant. 

The algorithm is designed with the purpose to make the states asymptotically converge to the consensus set of points, in which all points are equal. We say that consensus is reached (in the states) when this occurs. In the multi-agent systems setting where each agent $i$ has a local coordinate system given by orthogonal transformation of a global coordinate system, we may express Algorithm~\ref{eq:main_algo} using only local and relative information. 

\subsection{Compact representation of consensus algorithm}
If we define $X(k) \in \mathbb{S}(n,d)$ as the matrix for which $[X(k)]_i = x_i^T(k)$ for all $i$ and $k$, we may equivalently express \eqref{eq:main_algo} as 
\begin{align}
\label{eq:mainAlg3}
X(k+1) & ~= F(X(k)) = D(AX(k))AX(k), \text{ or } \\
\label{eq:mainAlg4}
x(k+1) & ~= f(x(k)) = (D(AX(k))A\otimes I_d)x(k),
\end{align} 
%where $X(k) \in \mathbb{S}(n,d)$ %and $[X(k)]_i = x_i(k)^T$ 
where $(\mathbb{S}^{d-1})^{n} \ni x(k) = \text{vec}(X(k)^T)$, and for a matrix $Z$ with $n$ non-zero rows: 
\begin{align}\label{eq:D}
D(Z) = &\text{diag} ( \| [Z]_1 \|_2 ^{-1}, \| [Z]_2 \|_2 ^{-1}, \ldots, \| [Z]_n \|_2 ^{-1}).
\end{align}
Depending on the context, we might use either $X$ or $x = \text{vec}(X^T)$ to represent the state of the system.  

\subsection{Fixed points and consensus points}\label{sec:periodic}
For dimension $d$ and a diagonally dominant weight matrix $A$, we define the set of fixed points for Algorithm \eqref{eq:mainAlg4} as
\begin{equation}
    \mathcal{F} = \{x \in (\mathbb{S}^{d-1})^n: f(x) = x\}.
\end{equation}
Of particular interest in this paper is the subset of fixed points comprising the consensus points  
\begin{equation}
    \mathcal{C} = \{x \in (\mathbb{S}^{d-1})^n: \exists \bar x \in \mathbb{S}^{d-1} \text{ s.t. } x = \text{vec}(\boldsymbol{1}_n\bar x^T).
\end{equation}
One can verify that $\mathcal{C}\subset \mathcal{F}$.
Indeed, if $x\in \mathcal{C}$, then the right-hand side of \eqref{eq:main_algo} becomes
\begin{align}
\frac{\sum_{j = 1}^n a_{ij}\bar x}{\sum_{j = 1}^n a_{ij}\|\bar x\|_2}=\bar x,
\end{align}
which implies that $f(x)=x$ for any $x \in \mathcal{C}$. 

\subsection{Tangent space, differential and unstable fixed points}
The tangent space of $\mathbb{S}^{d-1}$ at a point $y_i \in \mathbb{S}^{d-1}$ is the set $\mathcal{T}_{\mathbb{S}^{d-1}}(y_i) = \{P_{y_i}v: v \in \mathbb{R}^d\}$, where $P_{y_i} = I_d - y_iy_i^T$ is the matrix representation for the linear operator for projection onto the tangent space $\mathcal{T}_{\mathbb{S}^{d-1}}(y_i)$. 

For a point $y \in (\mathbb{S}^{d-1})^n$, there exists $Y \in \mathbb{S}(n,d)$ such that $y = \text{vec}(Y^T)$. Define $y_i^T = [Y]_i$ (the $i$'th row of $Y$). 
The tangent space of $(\mathbb{S}^{d-1})^{n}$ at the point $y \in (\mathbb{S}^{d-1})^n$ is then defined as
\begin{equation}
    \mathcal{T}_{(\mathbb{S}^{d-1})^n}(y) = \{P_yv: v \in \mathbb{R}^{nd}\},
\end{equation}
 where $P_{y}$ is a block-diagonal matrix with $n$ diagonal blocks, where the $i$'th diagonal block is $P_{y_i}$.
 
The differential $\mathcal{D}f(x)$ is the linear operator $ \mathcal{T}_{(\mathbb{S}^{d-1})^n}(x) \rightarrow \mathcal{T}_{(\mathbb{S}^{d-1})^n}(f(x))$, defined as follows~\cite{absil2009optimization,lee2019first}. For a curve $\gamma(t) \in (\mathbb{S}^{d-1})^n$ with $\gamma(0) = v \in \mathcal{T}_{(\mathbb{S}^{d-1})^n}(x)$,  $\mathcal{D}f(x)v = \frac{f \circ \gamma}{dt}(0) \in \mathcal{T}_{(\mathbb{S}^{d-1})^n}(f(x))$. Let $\text{det}(\mathcal{D}f(x))$ be the determinant of the $(n(d-1)) \times (n(d-1))$ matrix representing $\mathcal{D}f(x)$ with respect to an arbitrary choice of orthonormal bases. 

We define the set of unstable fixed points
\begin{align}
\mathcal{U}=\left\{ x \in \mathcal{F}: \max_i |\lambda_i (\mathcal{D}f(x))| > 1 \right\},
\end{align}
i.e., the condition for a fixed point to be unstable is that the largest in magnitude eigenvalue of the operator $\mathcal{D}f(x)$ is strictly larger than 1. 

\section{Properties of the differential}\label{sec:properties_diff}
The projected Jacobian matrix $J(x) \in \mathbb{R}^{nd \times nd}$, or simply $J$ for brevity, of $f(x)$ is used to represent the differential, i.e., the linear operator $\mathcal{D}f(x)$. In what follows we provide the expression for $J$ and then,  
for orthonormal bases of the tangent spaces, we provide the matrix $M \in \mathbb{R}^{n(d-1) \times n(d-1)}$, which is a matrix representation of the differential with respect to the bases.

Now, for a vector $v \in \mathcal{T}_{(\mathbb{S}^{d-1})^n}(x)$, 
\begin{equation}
    \mathcal{D}f(x)v = Jv \in \mathcal{T}_{(\mathbb{S}^{d-1})^n}(f(x)).
\end{equation}
To obtain a simplified expression of the projected Jacobian matrix $J$, let
$y = f(x)$, where $x, y \in (\mathbb{S}^{d-1})^n$. So $y$ is a function of $x$. Let $X, Y \in \mathbb{S}(n,d)$ be such that $y = \text{vec}(Y^T)$, $x = \text{vec}(X^T)$. Let $y_i^T = [Y]_i$ and $x_i^T = [X]_i$ for all $i$. Let $y = [y_1^T, y_2^T, \ldots, y_n^T]^T$ and $x = [x_1^T, x_2^T, \ldots, x_n^T]^T$.  By definition, $y_i=\frac{([A]_i \otimes I_d)x}{\|([A]_i \otimes I_d)x\|_2}$. 

The projected Jacobian matrix $J$ is a square block matrix with $n^2$ matrix blocks. We write $J = [J_{ij}]$, where $J_{ij} \in \mathbb{R}^{d \times d}$ for each $i,j \in \{1,2, \ldots, n\}$. Let $f_i(x)$ be the right-hand side of \eqref{eq:main_algo} for each $i$. Each $J_{ij}$-matrix is obtained by computation of the Euclidean gradient of each element in $f_i(x)$ followed by projection onto $\mathcal{T}_{\mathbb{S}^{d-1}}(x_j)$ using $P_{x_j}$: 
\begin{align}
    \nonumber
    J_{ij} =&  \frac{a_{ij}}{\|([A]_i \otimes I_d)x\|_2}\left(I_d - \frac{([A]_i \otimes I_d)xx^T([A]_i^T \otimes I_d)}{\|([A]_i \otimes I_d)x\|_2^2}\right) \\
     \nonumber
     ~&\cdot  (I_d - x_jx_j^T) \\
     %=&\frac{a_{ij}}{\|([A]_i \otimes I_d)x\|_2}(I_d - y_iy_i^T)(I_d - x_jx_j^T)\nonumber\\
     =& \frac{a_{ij}}{\|([A]_i \otimes I_d)x\|_2}P_{y_i}P_{x_j}.
\end{align} 
All together the $J_{ij}$-matrices comprise the projected Jacobian matrix 
\begin{equation}\label{eq:alf:1}
    J = P_y(D(AX)A \otimes I_d)P_x. 
\end{equation}
From this projected Jacobian matrix $J = [J_{ij}]$ (which depends on $x$), for chosen bases we can construct a matrix $M = [M_{ij}]$, which is the matrix representation of the differential $\mathcal{D}f(x)$ w.r.t. these bases. For each $i,j \in \{1,2,\ldots, n\}$,  $M_{ij} \in \mathbb{R}^{(d-1) \times (d-1)}$.  We proceed to obtain $M$ in the following manner. 

For each $i$, let $R_{y_i} \in \mathbb{R}^{d \times (d-1)}$ be a matrix whose columns are orthogonal to $y_i$ and which satisfies $R_{y_i}^TR_{y_i} = I_{d-1}$. The columns of $R_{y_i}$ comprise a basis for $\mathcal{T}_{\mathbb{S}^{d-1}}(y_i)$. Also, we let $R_{x_i} \in \mathbb{R}^{d \times (d-1)}$ be a matrix whose columns are orthogonal to $x_i$ and which satisfies $R_{x_i}^TR_{x_i} = I_{d-1}$. The columns of $R_{x_i}$ comprise a basis for $\mathcal{T}_{\mathbb{S}^{d-1}}(x_i)$. 
It holds that 
\begin{align}
    \label{eq:alf:10}
    P_{y_i} & = I_d - y_iy_i^T =  R_{y_i}R_{y_i}^T, \\
    \label{eq:alf:11}
    P_{x_i} & = I_d - x_ix_i^T =  R_{x_i}R_{x_i}^T.
\end{align}
Thus, we may write 
\begin{equation}
    J_{ij} = \frac{a_{ij}}{\|([A]_i \otimes I_d)x\|_2}R_{y_i}R_{y_i}^TR_{x_j}R_{x_j}^T.
\end{equation} 
Let $v_j = R_{x_j}\tilde{v}_j$, where $\tilde{v}_j \in \mathbb{R}^{d-1}$ is the $j$'th sub-vector of a total coordinate vector with respect to the basis for $\mathcal{T}_{(\mathbb{S}^{d-1})^n}(x)$. Let $u_{ij}=J_{ij}v_j$ and let $\tilde u_{ij} = R^T_{y_i}{u}_{ij} \in \mathbb{R}^{d-1}$ be the contribution of $\tilde v_{j}$ to $i$'th sub-vector of the coordinate vector for the basis chosen for $\mathcal{T}_{(\mathbb{S}^{d-1})^n}(y)$ when coordinates are mapped with the differential. It follows that
\begin{equation}
    \tilde{u}_{ij} = R_{y_i}^TJ_{ij}R_{x_j}\tilde v_{j} = \underbrace{\frac{a_{ij}}{\|([A]_i \otimes I_d)x\|_2}R_{y_i}^TR_{x_j}}_{ M_{ij}}\tilde v_{j} \in \mathbb{R}^{d-1}.
\end{equation}
Define $M = [M_{ij}]$. 
We may introduce the matrices $R_x \in \mathbb{R}^{nd \times n(d-1)}$ and $R_y  \in \mathbb{R}^{nd \times n(d-1)}$, where $R_x$ is the block-diagonal matrix in-which $R_{x_i}$ is the $i$'th block for $i = \{1,2, \ldots, n\}$, and $R_y$ is the block-diagonal matrix in-which $R_{y_i}$ is the $i$'th block for $i = \{1,2, \ldots, n\}$. The columns of $R_y$ comprise the basis for $\mathcal{T}_{(\mathbb{S}^{d-1})^n}(y)$, and the columns of $R_x$ comprise the basis for $\mathcal{T}_{(\mathbb{S}^{d-1})^n}(x)$. 

Given this new notation, we may rewrite \eqref{eq:alf:1} as
\begin{align}
\label{eq:alf:2}
    J & = R_yR_y^T(D(AX)A \otimes I_d)R_xR_x^T, \\
    \label{eq:alf:3}
    M & = R_y^T(D(AX)A \otimes I_d)R_x \\
    & = (D(AX) \otimes I_{d-1})R_y^T(A \otimes I_d)R_x, \nonumber
\end{align}
where $M$ is the matrix-representation of the differential with respect to the chosen bases. 

For any other orthogonal basis for $\mathcal{T}_{(\mathbb{S}^{d-1})^n}(y)$ there is a $Q_{y} \in \mathbb{O}(n(d-1))$ (orthogonal matrix) such that the columns of $\hat R_y = R_yQ_y$ comprise the basis, and for each orthogonal basis for $\mathcal{T}_{(\mathbb{S}^{d-1})^n}(x)$ there is a $Q_{x} \in \mathbb{O}(n(d-1))$ (orthogonal matrix) such that the columns of $\hat R_x = R_xQ_x$ comprise that basis. It can be seen that 
\begin{align}
    J & = R_yR_y^T(D(AX)A  \otimes I_d)R_xR_x^T \\ \nonumber
    & = R_y\underbrace{Q_yQ_y^T}_{I_{d-1}}R_y^T(D(AX)A  \otimes I_d)R_x\underbrace{Q_xQ_x^T}_{I_{d-1}}R_x^T\\ 
   &= \hat R_y\hat R_y^T(D(AX)A  \otimes I_d)\hat R_x \hat R_x^T. \nonumber
\end{align}
However, the matrix $M$, given by \eqref{eq:alf:3} for the old bases, will change to 
\begin{align}
    \nonumber
    \hat M  & = Q_y^TR_y^T(D(AX)A \otimes I_d)R_xQ_x \\
    \label{eq:alf:9}
    & = (D(AX) \otimes I_{d-1})Q_y^TR_y^T(A \otimes I_d)R_xQ_x 
\end{align}
for the new basis. But, as long as orientation is preserved by choosing $Q_x$ and $Q_y$ as matrices in $\mathbb{SO}(n(d-1))$, it holds that 
$\text{det}(\hat M) = \text{det}(Q_y^TR_y^T(D(AX)\!\otimes \! I_d)R_xQ_x) = \text{det}(R_y^T(D(AX)A \otimes I_d)R_x) = \text{det}(M)$. This entity is what is referred to as $\text{det}(\mathcal{D}f(x))$. 

Now we provide the first result of the paper, concerning $\text{det}(\mathcal{D}f(x))$. It is intuitive that making the diagonal elements of $A$ large enough ensures that the determinant is non-zero for all $x \in (\mathbb{S}^{d-1})^n$. However, the following proposition shows that the ``large enough'' is quite small: it suffices to increase the diagonal entries by a factor of only $\sqrt{2}$ beyond what is needed for strict diagonal dominance, independently of $n$ or $d$.  
 
\begin{proposition}\label{prop:nisse:2}
If $d \geq 2$ and $A$ is a weight matrix for a strongly connected directed graph $\mathcal{G}(n) = (\mathcal{V}(n), \mathcal{E})$
which satisfies 
    \begin{align}
       a_{ii} & >  \sqrt{2}\sum_{j \neq i}a_{ij}~\forall i,
    \end{align}
    then $\text{det}(\mathcal{D}f(x)) \neq 0$ for all $x\in (\mathbb{S}^{d-1})^n$.
\end{proposition}

\textbf{Proof:} We assume throughout that $A$ is strictly diagonally dominant. We begin by inspecting \eqref{eq:alf:3} and we note that $\text{det}(M) \neq 0 \Longleftrightarrow \text{det}(\tilde M) \neq 0$, where 
\begin{equation}
    \tilde M = R_y^T(A \otimes I_d)R_x. 
\end{equation}
Thus we can limit our attention to the the matrix $\tilde{M}$ in the continuation of this proof, where $\tilde{M}_{ij} = a_{ij}R_{y_i}^TR_{x_j}$. 

Let $y = f(x)$, where $x, y \in (\mathbb{S}^{d-1})^n$. Let $X, Y \in \mathbb{S}(n,d)$ be such that $y = \text{vec}(Y^T)$, $x = \text{vec}(X^T)$. Let $y_i^T = [Y]_i$ and $x_i^T = [X]_i$ for all $i$. Let $y = [y_1^T, y_2^T, \ldots, y_n^T]^T$ and $x = [x_1^T, x_2^T, \ldots, x_n^T]^T$. 

Without loss of generality, we can assume that $a_{ii} = 1$ for all $i$, since the dynamics of our system does not change when multiplying $A$ with a positive diagonal matrix from the left. To make this clear, we observe that 
\begin{align}
x_i(k+1) = \frac{\sum_{j = 1}^n a_{ij}x_j(k)}{\|\sum_{j = 1}^n a_{ij}x_j(k)\|_2} = \frac{\sum_{j = 1}^n \alpha a_{ij}x_j(k)}{\|\sum_{j = 1}^n \alpha a_{ij}x_j(k)\|_2}
\end{align}
for any $\alpha > 0$. Also, since $\text{det}(\mathcal{D}f(x))$ is invariant (does not change) under orientation-preserving change of orthonormal bases, we can select the $R_{x_i}$ matrices and the $R_{y_i}$ matrices, whose columns all together form the basis for the $\mathcal{T}_{(\mathbb{S}^{d-1})^n}(x)$ tangent spaces and the $\mathcal{T}_{(\mathbb{S}^{d-1})^n}(y)$ tangent spaces, respectively, in a certain way. We do so as follows. 

For each $i$, let $[R_{x_i}]^k$ and $[R_{y_i}]^k$ be the $k$'th column of $R_{x_i}$ and $R_{y_i}$, respectively, where $k \in \{1,2, \ldots, d-1$\}. For $d \geq 3$, we can choose $d-2$ unit vectors that are orthonormal and orthogonal to both $x_i$ and $y_i$. We let those vectors be the last $d-2$ columns of $R_{x_i}$ and $R_{y_i}$. So, $[R_{x_i}]^k = [R_{y_i}]^k$, $([R_{x_i}]^k)^T[R_{x_i}]^k = 1$,  $([R_{x_i}]^k)^T[R_{x_i}]^l = 0$, and $([R_{x_i}]^k)^Tx_i=([R_{x_i}]^k)^Ty_i=0$ for $k, l \geq 2$ and $k \neq l$. For $d = 2$, these constructed matrices are empty. 

Now, for each $i$, let $\cos(\theta_{i}) = y_i^Tx_i$. For a strictly diagonally dominant matrix $A$, it holds that $\cos(\theta_{i}) >0$ since 
\begin{equation}
    y_i^Tx_i\!=\!\sum_{j = 1}^n\frac{(a_{ij}x_j)^T}{\|\sum_{j = 1}^n a_{ij}x_j\|_2}x_i\!\geq\!\frac{a_{ii}-\sum_{j \neq  i} a_{ij}}{\|\sum_{j = 1}^n a_{ij}x_j\|_2}\!>\!0. 
\end{equation}
There is a linear subspace of $\mathbb{R}^d$ of dimension not larger than $2$  that contains $x_i$ and $y_i$ and is orthogonal to the linear subspace with basis comprising the last $d-2$ columns of $R_{x_i}$ or $R_{y_i}$. Let the orthogonal vectors $b_1$ and $b_2$ span this linear subspace. It holds that
\begin{equation}\label{eq:alf:4}
    \begin{bmatrix}
    x_i & y_i
    \end{bmatrix}
    = 
    \begin{bmatrix}
        b_1 & b_2
    \end{bmatrix}
    \begin{bmatrix}
        c_{x_i}^1 & c_{y_i}^1 \\
        c_{x_i}^2 & c_{y_i}^2
    \end{bmatrix}.
\end{equation}
It holds that $\|[c_{x_i}^1, c_{x_i}^2]^T\|_2 = \|[c_{y_i}^1, c_{y_i}^2]^T\|_2 = 1$. Furthermore, $\cos(\theta_{i}) = y_i^Tx_i = c_{x_i}^1c_{y_i}^1 + c_{x_i}^2c_{y_i}^2$. We select
\begin{equation}\label{eq:alf:5}
    \begin{bmatrix}
    [R_{x_i}]^1 & [R_{y_i}]^1
    \end{bmatrix}
    = 
    \begin{bmatrix}
        b_1 & b_2
    \end{bmatrix}
    \begin{bmatrix}
        c_{x_i}^2 & c_{y_i}^2 \\
        -c_{x_i}^1 & -c_{y_i}^1
    \end{bmatrix},
\end{equation}
whereby $([R_{x_i}]^1)^Tx_i = 0$, $([R_{y_i}]^1)^Ty_i = 0$, and $([R_{x_i}]^1)^T[R_{y_i}]^1=\cos(\theta_i)$.
All-in-all this means that we have selected the $R_{x_i}$-matrices and the $R_{y_i}$-matrices such that 
$\tilde M_{ii} = R_{y_i}^TR_{x_i} = \text{diag}([\cos(\theta_{i}), \boldsymbol{1}_{d-2}^T]^T)$.  

In what follows we make use of an argument that shares resemblance with the Gershgorin circle theorem, but here we consider blocks instead of rows of $\tilde M$. 
Let $\tilde M_i = [\tilde M_{i1}, \tilde M_{i2}, \ldots, \tilde M_{in}]$ for all $i$, and $\tilde M = [\tilde M_1^T, \tilde M_2^T, \ldots \tilde M_n^T]^T$. Thus, each $\tilde M_i$ contains $d-1$ rows of $\tilde M$. Let $v = [v_1^T, v_2^T, \ldots, v_n^T]^T$ be an eigenvector for an eigenvalue $\lambda \in \mathbb{C}$ of $ \tilde M$, where $v_i \in \mathbb{C}^{d-1}$ for all $i$. We normalize $v$ such that $\|v_{i_{\max}}\|_2 = 1$, where $i_{\max}=\arg\max\limits_{i}\{\|v_i\|_2\}$. 

Recall that we have chosen $a_{ii} = 1$ for all $i$. For the eigenvalue $\lambda$ with  eigenvector $v$, let $i = i_{\max}$.
We have that 
\begin{align}
    \lambda v_{i}  & = \tilde M_{ii}v_{i} + \sum_{j \neq i}\tilde M_{ij}v_{j}, \text{ and } \\
    \label{eq:alf:6}
    \text{Re}(\lambda)  & =\text{Re}(v_i^H \tilde M_{ii}v_{i}) + \sum_{j \neq i}\text{Re}(v_i^H\tilde M_{ij}v_{j}),
\end{align}
where $v_i^H$ denotes the conjugate transpose of $v_i$. By using the fact that $\tilde M_{ii} = \text{diag}([\cos(\theta_{i}), \boldsymbol{1}_{d-2}^T]^T)$, we conclude that $v_i^H\tilde M_{ii}v_{i}$ is real and greater or equal to $\cos(\theta_i)$. Furthermore, $|v_i^H\tilde M_{ij}v_{j}| = a_{ij}|v_i^HR_{y_i}^TR_{x_j}v_{j}|$ where $ |v_i^HR_{y_i}^TR_{x_j}v_{j}| \leq \|R_{y_i}v_i\|_2\|R_{x_i}v_j\|_2 = \|v_i\|_2\|v_j\|_2\leq \|v_i\|_2^2 \leq 1$. This means that $\text{Re}(v_i^HR_{y_i}^TR_{x_j}v_{j}) \geq -1$. By using these results together with \eqref{eq:alf:6}, we conclude that 
\begin{align}
    \text{Re}(\lambda) & \geq \cos(\theta_{i}) - \sum_{j \neq i}a_{ij}.
\end{align}
So, if 
\begin{equation} \label{eq:cos_ai}
   \cos(\theta_i) > \sum_{j \neq i}a_{ij} 
\end{equation}
is fulfilled for all $i$ and all $x \in (\mathbb{S}^{d-1})^n$, then 
$\text{det}(\mathcal{D}f(x)) \neq 0$ for all $x \in (\mathbb{S}^{d-1})^n$. 

Define $a_i = \sum_{j \neq i}a_{ij}$, which is strictly less than $1$ since $A$ is strictly diagonally dominant. It holds that 
\begin{align}
    \cos(\theta_i) & = x_i^Ty_i = x_i^T\frac{x_i +  \sum_{j \neq i}a_{ij} x_j}{\|x_i + \sum_{j \neq i} a_{ij}x_j\|_2}.
\end{align}
Due to the triangle inequality, $\|\sum_{j \neq i} a_{ij}x_j\|_2 \leq a_i$. Hence, there exist $\alpha_i \leq a_i$ and $z_i \in \mathbb{S}^{d-1}$ such that $\alpha_iz_i = \sum_{j \neq i} a_{ij}x_j$. Thus,
\begin{align}
    \cos(\theta_i) & = x_i^Ty_i \geq \min_{z \in \mathbb{S}^{d-1}}x_i^T\frac{x_i + \alpha_i z}{\|x_i + \alpha_i z\|_2} \\
    & = \min_{\phi_i}\frac{1 + \alpha_i\cos(\phi_i)}{(1 + 2\alpha_i\cos(\phi_i) + \alpha_i^2)^{\frac{1}{2}}},
\end{align}
where $\cos(\phi_i) = z^Tx_i$. Let
\begin{align}
    g(\phi_i) & = \frac{1 + \alpha_i\cos(\phi_i)}{(1 + 2\alpha_i\cos(\phi_i) + \alpha_i^2)^{\frac{1}{2}}}, \text{ whereby } \\
    \label{eq:alf:8}
    \frac{dg(\phi_i)}{d\phi_i} &  = -\frac{\alpha_i^2 \sin(\phi_i) (\alpha_i + \cos(\phi_i))}{(1 + 2 \alpha_i \cos(\phi_i) + \alpha_i^2)^{\frac{3}{2}}}.
\end{align}
By inspecting \eqref{eq:alf:8}, we obtain two candidates as optimal choices for $\phi_i$. For the first, $\phi_{\text{cand}, 1}$, it holds that $\cos(\phi_{\text{cand}, 1}) = - \alpha_i$, and for the second, $\phi_{\text{cand}, 2}$, it holds that $\sin(\phi_{\text{cand}, 2}) = 0$. For the first, it holds that 
$g(\phi_{\text{cand}, 1}) = \sqrt{1 - \alpha_i^2}$, and for the second it holds that $g(\phi_{\text{cand}, 2}) = 1$. So clearly $\phi_{\text{cand}, 1}$ is the minimizer. 

At this point we conclude that $\cos(\theta_i) \geq \sqrt{1 - \alpha_i^2} \geq \sqrt{1 - a_i^2}$. But we also had the requirement \eqref{eq:cos_ai} from before that $\cos(\theta_i) > a_i$. If we combine these two relations, we end up with the following sufficient condition for non-zero eigenvalues of $M$ 
\begin{equation} \label{eq:square2}
    \sqrt{1 - a_i^2} > a_i  \Longrightarrow    a_i < \frac{1}{\sqrt{2}}. 
\end{equation}
Thus, if \eqref{eq:square2} is satisfied for all $i \in \{1, \cdots, n\}$, then $\tilde M$ has no zero-eigenvalues, which means $M$ has no zero-eigenvalues.
\hfill $\qed$

We continue with the following result on the consensus set. 

\begin{proposition}\label{prop:nisse:1}
    For any strictly diagonally dominant weight matrix $A$ and $d \geq 2$,
    \begin{equation}
        \mathcal{C} \cap \mathcal{U} = \emptyset.
    \end{equation}
\end{proposition} 
\textbf{Proof:} See Appendix~\ref{proof:prop:nisse:1}. \hfill $\qed$

\section{Consensus for symmetric weight matrices} \label{sec:symmetric}
In this section we assume that the matrix $A$ used for Algorithm~\eqref{eq:main_algo} (equivalently expressed in \eqref{eq:mainAlg3} and \eqref{eq:mainAlg4}) is a symmetric strictly diagonally dominant weight matrix for a symmetric strongly connected graph $\mathcal{G}(n)$ (equivalent to connected undirected graph).

\subsection{Gradient ascent}
The following problems are equivalent: 
\begin{align}
(P_1)  \quad   & \min_{X \in \mathbb{S}(n,d)}\sum_{i,j}a_{ij}\|[X]_i^T - [X]_j^T\|_2^2, \\
 (P_2) \quad    & \max_{X \in \mathbb{S} (n,d)}\text{tr}(X^TAX), \\
 (P_3) \quad    & \max_{x \in (\mathbb{S}^{d-1})^n}x^T(A \otimes I_d)x. 
\end{align}
This is true since $[X]_i[X]_i^T = 1$ for all $i$ if $X \in \mathbb{S}(n,d)$, and $\text{tr}(X^TAX) = \sum_{i,j}a_{ij}[X]_i[X]_j^T$. Furthermore, $(P_3)$ is simply $(P_2)$ expressed using $x = \text{vec}(X^T)$. Our algorithm can be seen as projected descent for problem $(P_1)$ or projected gradient ascent with ``infinite step size'' for problems $(P_2)$ and $(P_3)$. 
    
We define the function $V_A(X)$ for $X \in \mathbb{S}(n,d)$ as 
\begin{align}
    V_A(X) =~& \text{tr}(X^TAX) = \sum_{i,j}a_{ij}[X]_i[X]_j^T.
\end{align}

\begin{lem}\label{lemma:lemma1}
Let $A$ be a weight matrix for a strongly connected graph. Then the following holds:
\begin{equation}
         \textnormal{arg}\max_{x \in (\mathbb{S}^{d-1})^n}x^T(A \otimes I_d)x = \mathcal{C}, 
    \end{equation}
{where $x = \text{vec}(X^T)$, $X \in \mathbb{S}(n,d)$.}
\end{lem}
\textbf{Proof:} See Appendix~\ref{proof:lemma:lemma1}. \hfill $\qed$

It holds that $F(X) = \text{arg}\max_{Y \in \mathbb{S}(n,d)}\text{tr}(Y^TAX)$.
Furthermore, if $A$ is a symmetric strictly diagonally dominant weight matrix (and thus positive definite), it follows that
\begin{align}
\nonumber
    0 \leq \text{tr}((F(X)-X)^TA(F(X)-X)) \\  
    \leq \text{tr}(F(X)^TAF(X))-\text{tr}(X^TAX),
\end{align}
where equality holds if and only if $\text{vec}(X^T)$ is a fixed point for Algorithm \eqref{eq:mainAlg4}, i.e., $X=F(X)$. Thus, if $X$ is not a fixed point, then $V_A(F(X)) > V_A(X)$. Since $V_A$ is also bounded over $\mathbb{S}(n,d)$, the sequence $\{V_A(F^k(X))\}_k$ converges.
Due to these facts, for $x \in (\mathbb{S}^{d-1})^n$, the set of limit points 
$\bigcap_{l = 0}^{\infty} \overline{\{f^k(x) : k \geq l\}} \subset \mathcal{F}$, i.e., the limit points are fixed points.

\subsection{Main results}
The following theorem, which is similar to the result for gradient flows in~\cite{johanl2017continious}, provides conditions on $A$ and $d$ such that all fixed points other than consensus points are unstable. 

\begin{proposition}\label{thm:sym:1}
If $d \geq 3$ and the $A$ that is used for Algorithm~\eqref{eq:main_algo} is a symmetric strictly diagonally dominant weight matrix for a symmetric connected graph $\mathcal{G}(n)$, then 
all fixed points not in the consensus set $\mathcal{C}$ are unstable. 
\end{proposition}

\textbf{Proof:} Let $x = \text{vec}(X^T)$ be a fixed point not in $\mathcal{C}$. We repeat in this proof that projected Jacobian matrix \eqref{eq:alf:1} at a fixed point is given by $J(x)=P_x(D(AX)A \otimes I_d)P_x=(D(AX)\otimes I_d)P_x(A \otimes I_d)P_x$, where the last equality holds since $P_x$ is a block-diagonal matrix.
Since $P_x^2=P_x$, the matrix $J(x)$ can be further rewritten as
\begin{equation}
    J(x)=(D(AX)\otimes I_d)H(x)+P_x,
\end{equation}
where $H(x)=P_x((A - D(AX)^{-1})\otimes I_d)P_x$ is symmetric. 

It holds that
\begin{align}
\nonumber
    & \exists i, \text{ s.t } \lambda_i (H(x)) > 0 \\ \Longrightarrow~&  \max_i |\lambda_i (J(x))| =  \max_i |\lambda_i (\mathcal{D}f(x))| > 1,
\end{align}
which means $x$ is an unstable fixed point. 

Indeed, assume that there exists a positive eigenvalue of the matrix $H(x)$. Then according to the Sylvester's law of inertia, there exists a positive eigenvalue for matrix $SH(x)S$, where $S=D(AX)^{1/2}\otimes I_d$. We denote this eigenvalue $\lambda_i$ and the corresponding right eigenvector $e_i$. It follows that $(D(AX)\otimes I_d)H(x)(Se_i)=(SSH(x)SS^{-1})(Se_i)=S(SH(x)S)e_i=\lambda_i(Se_i)$ and, thus, $Se_i$ is an eigenvector of $(D(AX)\otimes I_d)H(x)$ with the eigenvalue $\lambda_i>0$. Furthermore,  $P_xSe_i$ is also an eigenvector of $(D(AX)\otimes I_d)H(x)$ for the same eigenvalue $\lambda_i$ since $P_x^2=P_x$. Finally, $J(x)(P_xSe_i)=(\lambda_i+1)(P_xSe_i)$ meaning that the matrix $J(x)$ has a positive eigenvalue strictly larger than $1$. 

The matrix $H(x)$ is a special matrix in a class of Hessian matrices for which it is known that there is a positive eigenvalue~\cite{johanl2017continious}.
We provide the proof for sake of completeness. Let $z_d \in \mathbb{S}^{d-1}$, $z_{nd} = \boldsymbol{1}_n \otimes z_d \in \mathbb{R}^{nd}$. It holds that $z_{nd}^TH(x)z_{nd} = z_d^{T}\tilde{H}(x)z_d$, where $\tilde{H}(x)=(\bold{1}_n^T\otimes I_d)H(x)(\bold{1}_n\otimes I_d) \in \mathbb{R}^{d \times d}$. It holds
\begin{equation}\label{eq:nisse:100}
    \text{tr}(\tilde{H}(x)) = \sum_{i}\sum_{j \neq i}a_{ij}(d - 2 + \cos^2\theta_{ij}-(d-1)\cos\theta_{ij}),
\end{equation}
where $\cos(\theta_{ij}) = x_i^Tx_j$. This holds since $\text{tr}(P_{x_i}P_{x_i})=d-1$, and $\text{tr}(P_{x_i}P_{x_j})=d-2+\cos^2\theta_{ij}$ for $i \neq j$.

For $d > 2$ and at least one pair $(i,j)$ such that $\cos(\theta_{ij}) < 1$ (which must exist since $x \not \in \mathcal{C}$ and $A$ is a weight matrix for a strongly connected graph), $\text{tr}(\tilde{H}(x))>0$, and thus there is a positive eigenvalue for the symmetric matrix $\tilde H(x)$. But then we can choose $z_d$ such that $0 < z_d^{T}\tilde{H}(x)z_d = z_{nd}^TH(x)z_{nd}$. Consequently, we can conclude that the symmetric matrix $H(x)$ has a positive eigenvalue. 
\hfill $\qed$

The next theorem provides conditions under which the set of initial $x_0$'s for which convergence to a fixed point not in $\mathcal{C}$ occurs has measure zero. 

\begin{theorem}\label{thm:sym:2}
If $d \geq 3$ and the $A$ that is used for Algorithm~\eqref{eq:main_algo} is a weight matrix for a symmetric connected graph that satisfies $a_{ii} >  \sqrt{2}\sum_{j \neq i}a_{ij}~\forall i$,
 then the subset
 \begin{align} 
\{x \in (\mathbb{S}^{d-1})^n: \exists~x_{\infty} \in \mathcal{F}\backslash \mathcal{C} \text{ s.t. } \lim\limits_{k \rightarrow  +\infty}f^k(x) = x_{\infty}\} \end{align} of $(\mathbb{S}^{d-1})^n$ has measure zero.  
\end{theorem}

\textbf{Proof:}
Due to Proposition~\ref{prop:nisse:1} and  Proposition~\ref{thm:sym:1}, $\mathcal{F}\backslash \mathcal{U}= \mathcal{C}$. 
Furthermore, due to Proposition~\ref{prop:nisse:2}, $\text{det}(\mathcal{D}f(x)) \neq 0$ for all $x\in (\mathbb{S}^{d-1})^n$. The desired result now follows by application of Theorem~2 in \cite{lee2019first}.  \hfill $\qed$

\section{General weight matrices} \label{sec:general_weight_matricec}
%\subsection{A negative result for the unit circle}
If $d\geq 3$, Proposition~\ref{thm:sym:1} ensures that all fixed points not in $\mathcal{C}$ are unstable. However, the situation for $d=2$ is not mentioned. The following negative result provides a partial answer in terms of a condition for general weight matrices (not necessarily symmetric) ensuring there are fixed points not in $\mathcal{C} \cup \mathcal{U}$. A strategy as in the proof of Theorem~\ref{thm:sym:2} does not ensure convergence to points in $\mathcal{C}$ for all but a measure zero set of initial points.

\begin{proposition}
    \label{thm:nisse:4}
    Suppose $d = 2$ and the $A$ that is used for Algorithm~\eqref{eq:main_algo} is a strictly diagonally dominant weight matrix for a strongly connected graph $\mathcal{G}(n)$. 
    \begin{enumerate}
        \item If $x = \text{vec}(X^T)$ is a fixed point such that $x_i^Tx_j > 0$ for all $i,j$ such that $a_{ij} > 0$, where $x_i^T = [X]_i$ for all $i$, then $x \not \in \mathcal{U}.$ 
        \item For symmetric graphs, there exist (symmetric) $A$-matrices (and $n$) for which there are fixed points not in $\mathcal{C} \cup \mathcal{U}$.
    \end{enumerate}
\end{proposition}
\textbf{Proof:} See Appendix~\ref{proof:thm:nisse:4}. \hfill $\qed$

\subsection{Main result} \label{Sec:theorem2}
In the previous section it was concluded that for the unit circle (i.e., $d = 2$), there might be fixed points in neither $\mathcal{C}$ nor $\mathcal{U}$. In this section, we flip the perspective. Instead of investigating the properties of fixed points not in $\mathcal{C}$,  we ask ourselves when such fixed points exist in general. For a randomly chosen strictly diagonally dominant weight matrix for a given strongly connected graph $\mathcal{G}(n)$, what kind of fixed points exist in general? Theorem~\ref{th:prop_measure0} below states that for the unit circle ($d = 2$), only for a measure zero subset of strictly diagonally dominant weight matrices, Algorithm~\ref{eq:main_algo} has fixed points that are neither antipodal nor in $\mathcal{C}$.

We begin by defining the set  
$\mathcal{A}(\mathcal{G}(n))$ as the set of strictly diagonally dominant weight matrices for the strongly connected directed graph $\mathcal{G}(n)$.  Next, we define the set 
\begin{align}
\tilde{\mathcal{C}} = \{x \in  (\mathbb{S}^{d-1})^n: x = \text{vec}(X^T), \text{rank}(X) = 1\}.
\end{align}
If $x = \text{vec}(X^T) \in \tilde{\mathcal{C}}$, then for any $x_i = [X]_i$ and $x_j = [X]_j$, they are either equal or antipodal. It follows that $\mathcal{C} \subset \tilde{\mathcal{C}}$. Points in $\tilde{\mathcal{C}}\backslash \mathcal{C}$ are unstable fixed points.  We further use the notation $\mathcal{F}(A)$ to explicitly denote the set of fixed points for our algorithm (expressed by \eqref{eq:mainAlg4}) when the weight matrix $A$ is used.

\begin{theorem} \label{th:prop_measure0}
    If at least one of the two conditions below is satisfied, the subset of $A$-matrices in $\mathcal{A}(\mathcal{G}(n))$ 
    that satisfy
    $\mathcal{F}(A) \neq \tilde{\mathcal{C}}$, has measure zero.    \begin{enumerate}
        \item $d = 2$ and $\mathcal{G}(n)$ is a strongly connected directed graph. \\
        \item $ d \geq 2$ and $\mathcal{G}(n)$ is the complete graph. 
    \end{enumerate}
\end{theorem}
\textbf{Proof:} See Appendix~\ref{proof:th:prop_measure0}. \hfill $\qed$

%We end this section with an important remark on the proof of Theorem~\ref{th:prop_measure0}, which concerns symmetric matrices. 
Theorem~\ref{th:prop_measure0} states that if i) $d=2$ and $\mathcal{G}(n)$ is a strongly connected graph or ii) $d \geq 2$ and $\mathcal{G}(n)$ is the complete graph, there is only a measure zero set of $A$-matrices for which there are fixed points not in $\tilde{\mathcal{C}}$. Under any of these conditions, suppose $A$ is a symmetric strictly diagonally dominant weight matrix for the symmetric strongly connected graph $\mathcal{G}(n)$. Suppose for this chosen $A$, there are fixed points not in $\tilde{\mathcal{C}}$. Then we may $\epsilon$-disturb $A$ with an arbitrarily small $\epsilon$. For a randomly chosen $\epsilon$ that respects the zero-structure of $\mathcal{G}(n)$, such that $A + \epsilon$ is a strictly diagonally dominant weight matrix for $\mathcal{G}(n)$. The probability of $A + \epsilon$ having fixed points not in $\tilde{\mathcal{C}}$ is $0$. 

However, what happens if we restrict $\epsilon$ so that $A + \epsilon$ is symmetric? We define $\mathcal{A}_{\text{sym}}(\mathcal{G}(n))$ as the set of symmetric strictly diagonally dominant weight matrices for the symmetric connected graph $\mathcal{G}(n)$. The set $\mathcal{A}_{\text{sym}}(\mathcal{G}(n))$ is clearly a subset of $\mathcal{A}(\mathcal{G}(n))$. 
Already for the complete graph, the strategy deployed to prove Theorem~\ref{th:prop_measure0} fails when $\mathcal{A}_{\text{sym}}(\mathcal{G}(n))$ is considered instead of $\mathcal{A}(\mathcal{G}(n))$.

Furthermore, when $\mathcal{G}(n)$ is not the complete graph, as shown in Section~\ref{sec:simulations}, for symmetric weight matrices we often find fixed points $x = \text{vec}(X^T)$  whose rank is greater than $1$  (which we never do for non-symmetrical ones). This makes it comfortable to conjecture that a version of Theorem~\ref{th:prop_measure0} does not hold when the $A$-matrices are required to be in $\mathcal{A}_{\text{sym}}(\mathcal{G}(n))$ instead of the larger set $\mathcal{A}(\mathcal{G}(n))$. 

\section{Simulations}\label{sec:simulations} 
For randomly chosen $n$, $d$, strongly connected graphs and strictly diagonally dominant weight matrices, a hundred million simulations were conducted for Algorithm~\ref{eq:main_algo}. In every case, convergence to consensus was observed.
However, demonstrating the theoretical results regarding fixed points not in $\mathcal{C}$ through simulations is more challenging.
The approach taken in this context is to modify Algorithm \ref{eq:main_algo} to achieve convergence to non-consensus fixed points. While this modified algorithm does not guarantee the identification of all possible non-consensus fixed points, it provides empirical support for our theoretical findings.

In Section \ref{sec:symmetric}, we have shown that for undirected graphs, Algorithm \ref{eq:main_algo} ensures that the function $V_A(X)$ is non-decreasing. Thus, an algorithm minimizing the bounded $V_A(X)$ (equivalently, maximizing $-V_A$) on $\mathbb{S}(n,d)$, will yield a fixed point from $\mathcal{F}\setminus\mathcal{C}$. 
Let us consider the matrix $M_A=\alpha_AI_n - A$, where $\alpha_A >\max_i\left (\sum_{j=1}^n a_{ij} \right )$. It follows that $M_A=M_A^T$ is strictly diagonal dominant, and hence positive definite. Thus, Algorithm \ref{eq:main_algo} is well-defined (no division by zero) for the matrix $M_A$ instead of $A$. Any fixed point of Algorithm \ref{eq:main_algo} with matrix $A$ is a fixed point of Algorithm \ref{eq:main_algo} with matrix $M_A$. Furthermore, $V_M(X)=\alpha_A n-V_A(X)$. 

\begin{table}[htbp]
\centering
\begin{tabular} {c|ccccc}
\toprule
$(n,d)/m$ & 1 & 2 & 3 & 4 \\
\midrule
(3,2) & 0.7571 & 0.2429 & 0 & 0 \\
(6,2) & 0.3528 & 0.6472 & 0 & 0 \\
(4,3) & 0.2953 & 0.7046 & 0.0001 & 0 \\
(6,3) & 0.494 & 0.9256 & 0.25 & 0 \\
(7,4) & 0.0192 & 0.8802 & 0.1006 & 0 \\
(8,5)& 0.055 & 0.7647 & 0.2298 & 0 \\
\bottomrule
\end{tabular}
\caption{For symmetric  graphs and symmetric weight matrices: distribution over $10^4$ simulations of rank $m$ for non-consensus fixed points for different choices of $(n,d)$.}
\label{tab:example}
\end{table}

We conducted numerous simulations for Algorithm \ref{eq:main_algo} with matrix $M_A$ instead of $A$. Despite being designed for symmetric weight matrices, when the algorithm with $M_A$ instead of $A$ was applied to non-symmetric weight matrices, only fixed points in the set $\tilde{\mathcal{C}}$ were obtained in simulations, giving support to the result in Theorem~\ref{th:prop_measure0}. The discussion at the end of Section \ref{Sec:theorem2} highlights that the strategy to prove Theorem~\ref{th:prop_measure0} does not work for the complete graph if $A$ is restricted to be symmetric. However, our simulations suggest that there exist points in the set $\mathcal{F}\setminus\tilde{\mathcal{C}}$ in general when the graph is strongly connected and the weight matrix is symmetric.  

In Table~\ref{tab:example}, we provide simulation results performed for random symmetric strictly diagonally dominant weight matrices $A$ corresponding to strongly connected graphs, various dimensions $n$ and $d$, and various initial conditions $X(0)\in \mathbb{S}(n,d)$. In each simulation, Algorithm \eqref{eq:main_algo} with $M_A$ converged to a non-consensus fixed point. The table shows, for different choices $(n,d)$, the distribution of the rank $m$ for those fixed points over $10^4$ simulations with random $A$-matrices. From these simulations, it is clear that for the symmetric weight matrices the set $\mathcal{F}\setminus\tilde{\mathcal{C}}$ is not empty in general.

\section{Conclusions}\label{sec:conclusions}
This paper investigates the convergence of projection-based iterative or discrete-time consensus algorithms on the Euclidean unit sphere. The update equation for the algorithm comprises projection (onto the unit sphere) of conical combinations of neighboring states, where neighborhoods are defined by a directed strongly  connected symmetric or non-symmetric graph. For symmetric graphs, the algorithm reduces to projected gradient ascent.
The first result, that parallels earlier known results for gradient flows, establishes that for unit sphere dimension greater than or equal to $2$, if the weight matrix is such that each diagonal element is more than $\sqrt{2}$ larger than the sum of the other elements in the corresponding row, then stable fixed points are consensus points, and the set of initial points for which the algorithm converges to a non-consensus fixed point has measure zero. The second result of the paper is  that for 1) directed strongly connected graphs and the unit circle or for 2) the complete graph and any unit sphere dimension, only for a measure zero set of weight matrices there are fixed points for the algorithm  which do not have consensus or antipodal configurations. 

\begin{ack}                             
This work was partially supported by the Wallenberg AI, Autonomous Systems and Software Program (WASP) funded by the Knut and Alice Wallenberg Foundation and partly supported by the Swedish Research
Council under Grant 2019-04769. 
\end{ack}
\bibliographystyle{plain}         
\bibliography{refs.bib} 

@article{johanl2017continious,
  title={Almost global consensus on the $ n $-sphere},
  author={Markdahl, Johan and Thunberg, Johan and Gon{\c{c}}alves, Jorge},
  journal={IEEE Transactions on Automatic Control},
  volume={63},
  number={6},
  pages={1664--1675},
  year={2017},
  publisher={IEEE}
}

@article{lee2019first,
  title={First-order methods almost always avoid strict saddle points},
  author={Lee, Jason D and Panageas, Ioannis and Piliouras, Georgios and Simchowitz, Max and Jordan, Michael I and Recht, Benjamin},
  journal={Mathematical programming},
  volume={176},
  pages={311--337},
  year={2019},
  publisher={Springer}
}

@article{schenato2011average,
  title={Average timesynch: A consensus-based protocol for clock synchronization in wireless sensor networks},
  author={Schenato, L. and Fiorentin, F.},
  journal={Automatica},
  volume={47},
  number={9},
  pages={1878--1886},
  year={2011},
  publisher={Elsevier}
}

@inproceedings{reynolds1987flocks,
  title={Flocks, herds and schools: A distributed behavioral model},
  author={Reynolds, Craig W},
  booktitle={Proceedings of the 14th annual conference on Computer graphics and interactive techniques},
  pages={25--34},
  year={1987}
}

@article{vicsek1995novel,
  title={Novel type of phase transition in a system of self-driven particles},
  author={Vicsek, T. and Czirok, A. and  Eshel, B-J. and Cohen, I.  and Shochet, O. },
  journal={Physical Review Letters},
  volume={75},
  number={6},
  pages={1226},
  year={1995},
  publisher={APS}
}

@article{olfati2007consensus,
  title={Consensus and cooperation in networked multi-agent systems},
  author={Olfati-Saber, R. and Fax, J.A. and Murray, R.M.},
  journal={Proceedings of the IEEE},
  volume={95},
  number={1},
  pages={215--233},
  year={2007},
  publisher={IEEE}
}

@inproceedings{ren2005survey,
  title={A survey of consensus problems in multi-agent coordination},
  author={Ren, W. and Beard, R. W. and Atkins, E.M.},
  booktitle={American Control Conference, 2005. Proceedings of the 2005},
  pages={1859--1864},
  year={2005},
  organization={IEEE}
}

@article{cao2013overview,
  title={An overview of recent progress in the study of distributed multi-agent coordination},
  author={Cao, Yongcan and Yu, Wenwu and Ren, Wei and Chen, Guanrong},
  journal={IEEE Transactions on Industrial informatics},
  volume={9},
  number={1},
  pages={427--438},
  year={2012},
  publisher={IEEE}
}

@article{savazzi2020federated,
  title={{Federated learning with cooperating devices: A consensus approach for massive IoT networks}},
  author={Savazzi, Stefano and Nicoli, Monica and Rampa, Vittorio},
  journal={IEEE Internet of Things Journal},
  volume={7},
  number={5},
  pages={4641--4654},
  year={2020},
  publisher={IEEE}
}

@article{dorfler2014synchronization,
  title={Synchronization in complex networks of phase oscillators: A survey},
  author={D{\"o}rfler, Florian and Bullo, Francesco},
  journal={Automatica},
  volume={50},
  number={6},
  pages={1539--1564},
  year={2014},
  publisher={Elsevier}
}

@article{sarlette2009synchronization,
  title={Synchronization on the circle},
  author={Sarlette, Alain and Sepulchre, Rodolphe},
  journal={arXiv preprint arXiv:0901.2408},
  year={2009}
}

@article{markdahl2020high,
  title={{High-dimensional Kuramoto models on Stiefel manifolds synchronize complex networks almost globally}},
  author={Markdahl, Johan and Thunberg, Johan and Gon{\c{c}}alves, Jorge},
  journal={Automatica},
  volume={113},
  pages={108736},
  year={2020},
  publisher={Elsevier}
}

@inproceedings{kuramoto1975self,
  title={Self-entrainment of a population of coupled non-linear oscillators},
  author={Kuramoto, Yoshiki},
  booktitle={International symposium on mathematical problems in theoretical physics},
  pages={420--422},
  year={1975},
  organization={Springer}
}

@article{geshkovski2023mathematical,
  title={A mathematical perspective on transformers},
  author={Geshkovski, Borjan and Letrouit, Cyril and Polyanskiy, Yury and Rigollet, Philippe},
  journal={arXiv preprint arXiv:2312.10794},
  year={2023}
}

@article{townsend2020dense,
  title={Dense networks that do not synchronize and sparse ones that do},
  author={Townsend, Alex and Stillman, Michael and Strogatz, Steven H},
  journal={Chaos: An Interdisciplinary Journal of Nonlinear Science},
  volume={30},
  number={8},
  year={2020},
  publisher={AIP Publishing}
}

@article{mcrae2024benign,
  title={Benign landscapes of low-dimensional relaxations for orthogonal synchronization on general graphs},
  author={McRae, Andrew D and Boumal, Nicolas},
  journal={SIAM Journal on Optimization},
  volume={34},
  number={2},
  pages={1427--1454},
  year={2024},
  publisher={SIAM}
}

@incollection{absil2009optimization,
  title={Optimization algorithms on matrix manifolds},
  author={Absil, P-A and Mahony, Robert and Sepulchre, Rodolphe},
  booktitle={Optimization Algorithms on Matrix Manifolds},
  year={2009},
  publisher={Princeton University Press}
}

@inproceedings{tron2012intrinsic,
  title={Intrinsic consensus on {SO}(3) with almost-global convergence},
  author={Tron, Roberto and Afsari, Bijan and Vidal, Ren{\'e}},
  booktitle={2012 IEEE 51st IEEE Conference on Decision and Control (CDC)},
  pages={2052--2058},
  year={2012}
}

@article{thunberg2018dynamic,
  title={{Dynamic controllers for column synchronization of rotation matrices: a QR-factorization approach}},
  author={Thunberg, Johan and Markdahl, Johan and Gon{\c{c}}alves, Jorge},
  journal={Automatica},
  volume={93},
  pages={20--25},
  year={2018},
  publisher={Elsevier}
}

@article{sepulchre2010consensus,
  title={Consensus on nonlinear spaces},
  author={Sepulchre, Rodolphe},
  journal={IFAC Proceedings Volumes},
  volume={43},
  number={14},
  pages={1029--1039},
  year={2010},
  publisher={Elsevier}
}

@article{canale2008almost,
  title={{Almost global synchronization of symmetric Kuramoto coupled oscillators}},
  author={Canale, Eduardo and Monz{\'o}n, Pablo},
  journal={Systems Structure and Control},
  volume={8},
  pages={167--190},
  year={2008},
  publisher={InTech}
}

@inproceedings{canale2010complexity,
  title={On the complexity of the classification of synchronizing graphs},
  author={Canale, Eduardo and Monz{\'o}n, Pablo and Robledo, Franco},
  booktitle={International Conference on Grid and Distributed Computing},
  pages={186--195},
  year={2010},
  organization={Springer}
}

@inproceedings{canale2010wheels,
  title={The wheels: an infinite family of bi-connected planar synchronizing graphs},
  author={Canale, Eduardo A and Monz{\'o}n, Pablo A and Robledo, Franco},
  booktitle={2010 5th IEEE Conference on Industrial Electronics and Applications},
  pages={2204--2209},
  year={2010},
  organization={IEEE}
}

@inproceedings{monzon2005global,
  title={{Global considerations on the Kuramoto model of sinusoidally coupled oscillators}},
  author={Monz{\'o}n, Pablo and Paganini, Fernando},
  booktitle={Proceedings of the 44th IEEE Conference on Decision and Control},
  pages={3923--3928},
  year={2005},
  organization={IEEE}
}

@phdthesis{sarlette2009geometry,
  title={Geometry and symmetries in coordination control},
  author={Sarlette, Alain},
  year={2009},
  publisher={Universite de Liege (Belgium)}
}

@article{taylor2012there,
  title={{There is no non-zero stable fixed point for dense networks in the homogeneous Kuramoto model}},
  author={Taylor, Richard},
  journal={Journal of Physics A: Mathematical and Theoretical},
  volume={45},
  number={5},
  pages={055102},
  year={2012},
  publisher={IOP Publishing}
}

\newpage~\newpage

\appendix
\section{Proofs}
\subsection{Proof of Proposition~\ref{prop:nisse:1}}\label{proof:prop:nisse:1}
\textbf{Proof:} 
For any strictly diagonally dominant matrix, $\mathcal{C} \subset \mathcal{F}$. If $x \in \mathcal{C}$, there is an $\bar x \in \mathbb{S}^{d-1}$ such that $x_i = \bar x$ for all $i$, where $x = [x_1^T, x_2^T, \ldots, x_n^T]^T$. Thus, $P_{x}=I_n \otimes P_{\bar x}$ is a block-diagonal matrix with the same block $P_{\bar x}=I_d - \bar x \bar x^T$ on the diagonal. It follows
\begin{align}
    J(x) &\!=\!P_x(D(AX)A \otimes I_d)P_x
     \!=\!(D(AX)A \otimes P_{\bar x}),
\end{align}
since $P_{\bar x}^2 = P_{\bar x}$. 
{Since $x \in \mathcal{C}$, it holds that $[D(AX)]_{ii}=1/\sum_{j=1}^n a_{ij}$. Thus,} $D(AX)A$ is a right-stochastic matrix, and by using the property of eigenvalues of the Kronecker product, we conclude that
\begin{align}
     & \max_i |\lambda_i(\mathcal{D}f(x))| = \max_i |\lambda_i(J(x))| \\
     =~& \max_i |\lambda_i(D(AX)A \otimes P_{\bar x})| = 1. \nonumber
\end{align}
\hfill $\qed$

\subsection{Proof of Lemma~\ref{lemma:lemma1}}\label{proof:lemma:lemma1}
\textbf{Proof:} 
Let $x = \text{vec}(X^T)$, then $-1 \leq [X]_i[X]_j \leq 1$. Suppose $x \neq \mathcal{C}$, then there is $(i,j) \in \mathcal{E}$ (i.e., the edge set of the connected graph) such that $[X]_i[X]_j < 1$ and $a_{ij} > 0$, whereby 
\begin{align}
    x^T(A \otimes I_d)x = \sum_{i,j}a_{ij}[X]_i[X]_j^T < \sum_{i,j}a_{ij}.
\end{align}
Suppose $x \in \mathcal{C}$, then 
\begin{align}
    x^T(A \otimes I_d)x = \sum_{i,j}a_{ij}[X]_i[X]_j^T = \sum_{i,j}a_{ij}.
\end{align}
\hfill $\qed$

\subsection{Proof of Proposition~\ref{thm:nisse:4}}\label{proof:thm:nisse:4}
\textbf{Proof:}
i) The matrix $M$, see \eqref{eq:alf:3}, representing $\mathcal{D}f(x)$ is now of size $n \times n$ and can be chosen in a certain way.  Let 
\begin{equation}
    R_{x_i} = \begin{bmatrix}
        0 & 1 \\
        -1 & 0
    \end{bmatrix}x_i,
\end{equation}
whereby 
\begin{equation}
    M_{ij} = \frac{a_{ij}R_{x_i}^TR_{x_j}}{\sum_{j}a_{ij}x_i^Tx_j} =\frac{a_{ij}x_i^Tx_j}{\sum_{j}a_{ij}x_i^Tx_j}~\forall i,j.
\end{equation}
  This is a right-stochastic matrix whose largest eigenvalue is $1$.
  However, the eigenvalues of $M = [M_{ij}]$ are the same as the nonzero eigenvalues of the Jacobian matrix $J = [R_{x_i}M_{ij}R_{x_j}^T]$. Indeed, since orthogonal similarity preserves eigenvalues, we can consider the eigenvalues of $Q^TJQ$ instead, where $Q = \text{diag}(Q_{x_1},Q_{x_2}, \ldots, Q_{x_n})$ is an orthogonal matrix with $Q_{x_i} = [R_{x_i},x_i]$.  It holds that
\begin{equation}
Q_{x_i}^TR_{x_i}M_{ij}R_{x_j}^TQ_{x_j} = \begin{bmatrix} 1 \\ 0\end{bmatrix}M_{ij} \begin{bmatrix} 1 & 0\end{bmatrix}=
M_{ij}\begin{bmatrix} 1 & 0\\
0 & 0\end{bmatrix},
  \end{equation}
and thus
\begin{align}
Q^TJQ = M \otimes \begin{bmatrix} 1 & 0\\
0 & 0\end{bmatrix}.  
\end{align}
Since the eigenvalues of the Kronecker product are given by all pairwise products of the eigenvalues of the factors, we conclude that $\max_i |\lambda_i (J)| = 1$.
 
ii) Let for example $n = 5$ and $A$ be such that $[A]_{ii} = 3, [A]_{i,\gamma(i-1)} = [A]_{i,\gamma(i+1)} = 1$ for all $i$, where $\gamma(0) = n, \gamma(n+1) = 1, \gamma(i) = i$~  for $i \in \{2,3, \ldots, n-1\}$. Let $x_0 = \text{vec}(X_0^T)$. Let $[X_0]_i = [\cos(2\pi(i-1)/5), \sin(2\pi(i-1)/5)]$ for all $i$. Then, $\text{vec}(X_0^T)$ is a fixed point for Algorithm \eqref{eq:mainAlg4} with the considered symmetric $A$. This constructed example satisfies condition i) and, thus, $\max_i |\lambda_i (\mathcal{D}f(x))| = 1$.
\hfill  $\qed$

We see that the example used in the proof of ii) in Proposition~\ref{thm:nisse:4} also satisfies i) in the proposition. A more trivial example that satisfies i) is when all the $x_i$'s are in the positive orthant whereby the fixed point $x$ must be in $\mathcal{C}$. 

\subsection{Proof of Theorem~\ref{th:prop_measure0} and a remarks about symmetric weight matrices}\label{proof:th:prop_measure0}
Before we provide the proof of Theorem~\ref{th:prop_measure0}, we introduce some notations. 
We assume throughout that $\mathcal{G}(n)$ is a strongly connected directed graph. 
We define $B(\mathcal{G}(n))$ as the unique binary matrix that acts as a weight matrix for $\mathcal{G}(n)$ and defines its zero-structure.  

For $\mathcal{S} \subset (\mathbb{S}^{d-1})^n$ (i.e., any subset of $(\mathbb{S}^{d-1})^n$ including $(\mathbb{S}^{d-1})^n$), we define 
\begin{align}   \mathcal{S}_{\text{prod}}(\mathcal{S}) & = \mathcal{A}(\mathcal{G}(n)) \times \{\text{diag}(v) : v \in (\mathbb{R}^+)^n \} \times \mathcal{S}.
\end{align}
Furthermore for $(A, D, x) \in \mathbb{R}^{n \times n} \times \{\text{diag}(v) : v \in (\mathbb{R}^+)^n \} \times \mathbb{R}^{n d}$ we define the function 
\begin{align} \label{eq:vectorized_for_jacobian_Ym}
    g(A,D,x) =  (A \otimes I_d) x-(D \otimes I_d )x.
\end{align}
We further introduce $\mathcal{F}(A,\mathcal{S}) \subset \mathcal{F}(A,(\mathbb{S}^{d-1})^n) = \mathcal{F}(A)$, which is the set of fixed points in $\mathcal{S} \subset (\mathbb{S}^{d-1})^n$ of Algorithm \eqref{eq:mainAlg4} for the weight matrix $A$. 

If $x \in \mathcal{F}(A, \mathcal{S})$, it is clear that $x$ is also in the set  
\begin{align}
  {\mathcal{F}_g}(A, \mathcal{S}) & {=} \{ x \in \mathcal{S}:\\ 
    &\exists~D \in  \{\text{diag}(v) : v \in (\mathbb{R}^+)^n \},
    & g(A,D,x) = 0\},\nonumber
\end{align}
i.e., $\mathcal{F}(A, \mathcal{S}) \subset \mathcal{F}_g(A, \mathcal{S})$.

We define
\begin{align}
\mathcal{S}_{\text{prod},g}(\mathcal{S}) & = \{(A,D,x) \in \mathcal{S}_{\text{prod}}(\mathcal{S}) : g(A,D,x) = 0\},
\end{align}
where-after we define 
\begin{align}   \mathcal{A}_g(\mathcal{G}(n), \mathcal{S}) =~ &  \{A \in \mathcal{A}(\mathcal{G}(n)) : \exists~D, x  \text{ s.t. } \\ \nonumber
    & (A,D,x) \in \mathcal{S}_{\text{prod},g}(\mathcal{S})\}. 
\end{align}
Given the above, if we prove for $d=2$ or for complete graph $\mathcal{G}(n)$ and $d\geq 2$, that the set $\mathcal{A}_g(\mathcal{G}(n), (\mathbb{S}^{d-1})^n \backslash \tilde{\mathcal{C})}$ has measure zero in $\mathcal{A}(\mathcal{G}(n))$, then we have proven Theorem~\ref{th:prop_measure0}.

\textbf{Proof of Theorem~\ref{th:prop_measure0}:} \\
To begin with, we note that for any strictly diagonally dominant $A$, $\tilde{\mathcal{C}} \subset \mathcal{F}(A)$. 
Second, for any subset $\mathcal{S}$ of $(\mathbb{S}^{d-1})^n$, the set $\mathcal{A}_g(\mathcal{G}(n), \mathcal{S})$ is not empty. 
This can be shown as follows. For $x = \text{vec}(X^T) \in \mathcal{S}$ and a positive matrix $D$, $g(A,D,x) = 0$ can be alternatively expressed as
\begin{equation}
    (I_n \otimes X^T)K_A\text{vec}(A^T) = (I_n \otimes X^T)K_D\text{vec}(D).
\end{equation}
The matrices $K_A = \text{diag}(\text{vec}(B(\mathcal{G}(n))^T))$ and $K_{D}=\text{diag}(\text{vec}(I_n))$ are $n^2 \times n^2$ diagonal (binary) projection matrices. 
However, $K_D = K_A \odot K_D$ (element-wise product), since $A$ is assumed diagonally dominant. Thus, there is always an $A$ that solves $g(A,D,x) = 0$ when $x \in \mathcal{S}$ and $D$ is a positive diagonal matrix. 

Now we continue by first addressing condition ii) and secondly addressing condition i). 

\textbf{ii) }{($\mathcal{G}(n)$ is the complete graph)} \\
We first define 
\begin{align}    
\mathcal{S}_m = \{x \in (\mathbb{S}^{d-1})^n: x = \text{vec}(X^T), \text{rank}(X) = m\}
\end{align}
for $m \in \{1,2, \ldots, \min\{d,n\}\}$. We note that $\mathcal{S}_1 = \tilde{\mathcal{C}}$, $(\mathbb{S}^{d-1})^n = \bigcup_{m = 1}^d\mathcal{S}_m$, and $(\mathbb{S}^{d-1})^n \backslash \tilde{\mathcal{C}}  = \bigcup_{m = 2}^d\mathcal{S}_m$. Furthermore, 
\begin{align}    
\mathcal{A}_g(\mathcal{G}(n), (\mathbb{S}^{d-1})^n \backslash \tilde{\mathcal{C}}) =  \bigcup_{m = 2}^d\mathcal{A}_g(\mathcal{G}(n), \mathcal{S}_m),
\end{align}
and for given $A$,
\begin{align}   
\mathcal{F}(A, (\mathbb{S}^{d-1})^n \backslash \tilde{\mathcal{C}}) \subset  \bigcup_{m = 2}^d\mathcal{F}_g(A,\mathcal{S}_m).\end{align}
 Now, for $x \in \mathcal{S}_m$, where $x = \text{vec}(X^T)$, there exists an orthogonal transformation matrix $R \in \mathbb{O}(d)$ such that all the elements of the last $d - m$ columns of $XR$ are zero and $[XR]_1 = [1, 0, 0, \ldots, 0]$. Furthermore, if $x \in \mathcal{F}(A, \mathcal{S}_m)$, so does $\text{vec}((XR)^T)$.  Hence, if we neglect last $d-m$ zero columns of $XR$ and define  
\begin{align}
\tilde{\mathcal{S}}_m = & \{x \in (\mathbb{S}^{m-1})^n: x = \text{vec}(X^T), \\
& [X]_1 = [1, 0, 0, \ldots, 0], \text{rank}(X) = m \}\nonumber
\end{align}
for $2 \leq m \leq \min\{d,n\}$, it holds that 
$\mathcal{A}_g(\mathcal{G}(n), \tilde{\mathcal{S}}_m)  = \mathcal{A}_g(\mathcal{G}(n),{\mathcal{S}}_m)$. 
Also,  $\tilde{\mathcal{S}}_m$ is a manifold of dimension $(n-1)(m-1)$. 

Let $2 \leq m \leq \min\{d,n\}$. The first order approximation of $g$ on $\mathcal{S}_{\text{prod}}(\tilde{\mathcal{S}}_m)$ at a point $(A, D, x) \in \mathcal{S}_{\text{prod}}(\tilde{\mathcal{S}}_m)$ is 
\begin{align}
    & g(A + \Delta_A,D + \Delta_D,x + \Delta_x) \\
    = & g(A,D,x) + J_{g,A}\text{vec}(\Delta_A^T) + J_{g,D}\text{diag}(\Delta_D) + J_{g,x}\Delta_x,\nonumber
\end{align}
where $J_{g,A}=(I_n \otimes X^T )K_A \in \mathbb{R}^{nm \times n^2}$, $J_{g,x}=((A-D) \otimes I_m)P_x \in \mathbb{R}^{nm \times nm}$, $J_{g,D}=-(I_n \otimes X^T )K_{D} \in \mathbb{R}^{nm \times n^2}$, and $[\Delta_x]_1 = \bold{0}_{m}^T$. 

We define $J_g(A,D,x) \in \mathbb{R}^{nm \times (nm+2n^2)}$ or, for notational convenience, simply $J_g$, as the matrix \begin{align}    
J_g = [J_{g,A}, J_{g,D}, J_{g,x}].
\end{align}
Since $\mathcal{G}(n)$ is assumed to be the complete graph, $B(\mathcal{G}(n))$ is the matrix whose all entries are equal to $1$. Thus, $K_A = I_{n^2}$ and $J_{g,A} = I_n \otimes X^T$. Since $X$ has rank $m$, this implies that $J_{g,A}$ has rank $mn$, and consequently $J_g$ has full rank. Thus,  according to the implicit function theorem, $\mathcal{S}_{\text{prod},g}(\tilde{\mathcal{S}}_m)$ is a manifold of dimension 
\begin{align}
    \nonumber
   \text{dim}(\mathcal{S}_{\text{prod},g}(\tilde{\mathcal{S}}_m)) =~ & n^2 + n + (n-1)(m-1) - nm \\
   =~ & n^2 - m + 1.
\end{align}
Let 
\begin{align}   
 h_m: \mathcal{S}_{\text{prod},g}(\tilde{\mathcal{S}}_m) \mapsto \mathcal{A}(\mathcal{G}(n))
 \end{align}
be defined as $h_m(A,D,x) = A$, and $h_m(\mathcal{S}_{\text{prod},g}(\tilde{\mathcal{S}}_m)) = \mathcal{A}_g(\mathcal{G}(n), \tilde{\mathcal{S}}_m)$. Since the dimension of the manifold $\mathcal{A}(\mathcal{G}(n))$ is $n^2 > n^2 - m + 1 = \text{dim}(\mathcal{S}_{\text{prod},g}(\tilde{\mathcal{S}}_m))$, the differential of $h_m$ has rank smaller than $n^2$ for all points in $\mathcal{S}_{\text{prod},g}(\tilde{\mathcal{S}}_m)$. Thus, {Sard's} theorem asserts that $\mathcal{A}_g(\mathcal{G}(n), \tilde{\mathcal{S}}_m)$ has measure zero in $\mathcal{A}_g(\mathcal{G}(n))$. 
This is true for each $m$ between $2$ and $d$. Hence, the finite union 
\begin{align}    
\bigcup_{m = 2}^d\mathcal{A}_g(\mathcal{G}(n), \tilde{\mathcal{S}}_m) = \mathcal{A}_g(\mathcal{G}(n), (\mathbb{S}^{d-1})^n \backslash \tilde{\mathcal{C}})
\end{align}
has measure zero in $\mathcal{A}(\mathcal{G}(n))$.

\textbf{i)} ($d = 2$ and $\mathcal{G}(n)$ is a strongly connected graph)\\
We borrow notation from the proof of ii). We want to prove that $\mathcal{A}_g(\mathcal{G}(n), \tilde{\mathcal{S}}_2)$ has measure zero in $\mathcal{A}(\mathcal{G}(n)).$

We may create a collection of sets $\{\tilde{\mathcal{S}}^k_2\}_{k}^{N_p}$, the union of which is $\mathcal{A}_g(\mathcal{G}(n), \tilde{\mathcal{S}}_2)$, i.e., $\bigcup_{k = 1}^{N_p}\tilde{\mathcal{S}}^k_2 = \mathcal{A}_g(\mathcal{G}(n), \tilde{\mathcal{S}}_2)$. 
This collection of sets is defined according to the procedure below. 

We first create all partitions of the set $\{1,2, \ldots, n\}$ except the trivial partition containing only $\{1,2, \ldots, n\}$. For each such partition, for example $\{\{1,5,3\}, \{2,4\}, \\ \{6\},\{7, 8, \ldots, n\}\}$, we create $2^n$ sets by allowing for change of sign of each integer, i.e., for the example above we may choose $\{\{1,-5,3\}, \{-2, 4\}, \{6\},\{-7, 8, \ldots, n\}\}.$ There are $2^n$ such sign-patterns for each partition.
Each $\tilde{\mathcal{S}}^k_2$ in the collection of sets corresponds to a partition and a choice of sign-pattern in the sense that all the $x_i$'s inside a set in the partition are equal up to a sign dictated by the pattern. Furthermore, as an additional constraint we require the following. If $x_i$ belongs to one set in the partition and $x_j$ belongs to another set in the partition, the matrix $[x_i, x_j]$ has full rank (i.e., the matrix has rank $2$, and $x_i$ and $x_j$ are linearly independent). 

To illustrate, for the partition and the sign pattern above, the $\tilde{\mathcal{S}}^k_2$-set is   
\begin{align}
    \{x \in (\mathbb{S}^1)^n: &  ~x = \text{vec}(X^T), \text{rank}(X) = 2,\\\nonumber
    & [X]_1 = -[X]_5 = [X]_3 = [1,0]\\\nonumber
   - & [X]_2 = [X]_4 \\ \nonumber
   - & [X]_7 = [X]_8 = \cdots = [X]_n \},
\end{align}
where furthermore each vector in one set in the partition is linearly independent of any other vector in another set of the partition. 
In this example the partition comprised four sets, and the $\tilde{\mathcal{S}}^k_2$-set comprises a manifold of dimension $3$. To get this dimension, we note that all the vectors in the set in the partition containing $x_1$ are equal to $[1, 0]^T$ up to sign, which adds $0$ to the total dimension. The remaining three sets in the partition each adds $1$ to the total dimension. 
In general, if the partition comprises $N^k_s$ sets, the $\tilde{\mathcal{S}}^k_2$-set is a manifold of dimension $N^k_s - 1$. Finally, the reason for excluding the trivial partition with only one set $\{1,2, \ldots, n\}$ together with all its sign-patterns, is that this partition with the sign-patterns corresponds to the set $\tilde{\mathcal{C}}$.

Now, let $\mathcal{N}_i = \{j: [B(\mathcal{G}(n))]_{ij} = 1\}$. This set is usually, in the multi-agent systems context, referred to as the set of neighbors of agent $i$. It comprises $i$ and the neighbors of $i$ in the graph $\mathcal{G}(n)$.  Suppose there is an $i$ such that $\mathcal{N}_i$ is a subset of a set in the partition $\tilde{\mathcal{S}}^k_2$. Then there is a vector $\bar x = [\cos(\theta), \sin(\theta)]^T \in \mathbb{S}^1$, such that $x_j = s_j\bar x $, where $s_j \in \{-1, 1\}$, which comes from the sign-pattern for $\tilde{\mathcal{S}}^k_2$, and $\theta \in [0, 2\pi)$ for all $j \in \mathcal{N}_i$. The two scalar equations represented by $[A]_iX -d_i[X]_i = 0$, where $d_i = [D]_{ii}$, can be written as 
\begin{align}
    &\cos(\theta)(\sum_{j\in \mathcal{N}_i}a_{ij}s_j - d_i) = 0,\\
    &\sin(\theta)(\sum_{j\in \mathcal{N}_i}a_{ij}s_j - d_i) = 0.
\end{align}
We see that we can replace these two equations by one equation
\begin{equation}
    \sum_{j\in \mathcal{N}_i}a_{ij}s_j - d_i = 0.
\end{equation}
So one of the previous two equations were redundant. 

In general, for each $\tilde{\mathcal{S}}^k_2$-set, there is a maximum number of equations $N^k_e$ 
of the $2n$ number of (scalar) equations in $g(A,D,x) = 0$ that are redundant. Removal of these equations makes the matrix with $2n - N^k_e$ rows, corresponding to $J_{g,A}$ in the proof of ii), full rank. Now, since the graph $\mathcal{G}(n) = (\mathcal{V}(n), \mathcal{E})$ is strongly connected, there is a loop path/sequence $\{i_l\}_{l = 1}^{n+1}$, where $i_1 = i_n+1$ and $\cup_{l = 1}^n\{i_l\}  = \mathcal{V}(n)$, such that $[B(\mathcal{G}(n))]_{i_{l}i_{l+1}} = 1$ for $l \in \mathcal{V}(n)$. We recall that $\mathcal{V}(n) = \{1,2, \ldots, n\}$. This means that for each set in the partition, there must be a node in the graph with a neighbor that is not in the node's set of the partition. For such nodes there are no redundant equations.
Thus, the reduction of the number of equations represented by the number $N_e$ satisfies 
\begin{equation}
    N^k_e \leq n - N^k_s.
\end{equation}
Now, according to the implicit function theorem,  $\mathcal{S}_{\text{prod},g}(\tilde{\mathcal{S}}^k_2)$ is a manifold of dimension 
\begin{align}
    \nonumber
    & \text{dim}(\mathcal{S}_{\text{prod},g}(\tilde{\mathcal{S}}^k_2)) \\
    \nonumber
    =~& \sum_{i,j}[B(\mathcal{G})]_{ij} +n + (N^k_s - 1)  - (2n-N^k_e) \\
    \nonumber
    \leq~& \sum_{i,j}[B(\mathcal{G})]_{ij} +n + (N^k_s - 1)  - (2n-(n - N_s^K)) \\
    =~&    \sum_{i,j}[B(\mathcal{G})]_{ij} -1,
\end{align}
 whereas
$\mathcal{A}(\mathcal{G}(n))$ is a manifold of dimension $\sum_{i,j}[B(\mathcal{G})]_{ij}$. Thus we may use Sard's theorem analogously as for the the last part of the proof of ii) to conclude that $\mathcal{A}_g(\mathcal{G}(n), \tilde{\mathcal{S}}^k_2)$ has measure zero in $\mathcal{A}(\mathcal{G}(n))$. Hence,
\begin{align}    
\mathcal{A}_g(\mathcal{G}(n), (\mathbb{S}^{1})^n \backslash \tilde{\mathcal{C}}) = \bigcup_{k = 1}^{N_p}\mathcal{A}_g(\mathcal{G}(n), \tilde{\mathcal{S}}^k_2)\end{align}
has measure zero in $\mathcal{A}(\mathcal{G}(n)).$
\hfill $\qed$

Now we show that, under the restriction that $\mathcal{G}(n)$ is the complete graph, the strategy deployed to prove Theorem~\ref{th:prop_measure0} fails when $\mathcal{A}_{\text{sym}}(\mathcal{G}(n))$ is considered instead of $\mathcal{A}(\mathcal{G}(n))$. Since for the complete graph and fixed point $x = \text{vec}(X^T)$, the matrix $J_g$ in the proof of Theorem~\ref{th:prop_measure0} has rank strictly less than $nm$ if $m>1$, where $m=\text{rank}(X)$. This can be shown as follows.  

 For a symmetric matrix $A$, the corresponding part of the matrix $J_g$ must be projected onto the tangent space of symmetric matrices, yielding
 \begin{align}
    J_{g,A}=(I_n \otimes X^T)D_n, 
 \end{align} 
 where $D_n \in \mathbb{R}^{n^2 \times \frac{n(n+1)}{2}}$ is a duplication matrix which maps unique elements of $A$ to its full vectorized form. The duplication matrix $D_n$ for a symmetric matrix $C$ is defined such that
\begin{align} \label{eq:dupl}
    D_n\textnormal{vech}(C)=\textnormal{vec}(C),
\end{align} where $\textnormal{vec}(C)\in \mathbb{R}^{n^2 \times 1}$, and $\textnormal{vech}(C)\in \mathbb{R}^{\frac{n(n+1)}{2} \times 1}$ stacks only the lower-triangular part (including the diagonal) of $C$ into a vector.
 
Let consider vector $w_0 \in \mathbb{R}^{n^2 \times 1}$ such that 
\begin{align}
    & w_0=\textnormal{vec}(B), \text{and}\\  
    & D_n^Tw_0=0, \label{eq:Dnwo}
\end{align}
where $B $ - some matrix in $\mathbb{R}^{n \times n}$. Then due to \eqref{eq:dupl}, for any symmetric matrix $C \in \mathbb{R}^{n \times n}$, it holds:
\begin{align}
0&=\textnormal{vech}^T(C)D^T_nw_0=\textnormal{vec}^T(C)\textnormal{vec}(B)=\textnormal{trace}(C^TB)\\ \nonumber
&=\textnormal{trace}(B^TC)=\textnormal{trace}(BC).
\end{align}
Thus, $\textnormal{trace}((B^T+B)C)=0$ for any symmetric matrix $C \in \mathbb{R}^{n \times n}$, implying that $B$ must be skew-symmetric in order to satisfy \eqref{eq:Dnwo}. 
Let consider skew-symmetric $B$ in the form $B=XR_0X^T$, where $R_0 \in \mathbb{R}^{m \times m}$ is skew-symmetric, and $X \in \mathbb{S}(n,m)$ such that $\text{vec}(X^T)$ is a fixed point. Then, it follows 
\begin{align}
    0=D_n^T\text{vec}(XR_0X^T)=D_n^T(I_n\otimes X)\text{vec}(R_0X^T),
\end{align} meaning that $\textnormal{vec}^T(R_0X^T)$ is a left eigenvector for $J_{g,A}J_{g,A}^T$ with corresponding zero eigenvalue. Overall, there are $\frac{m(m-1)}{2}$ independent skew-symmetric matrices $R_0$ defining $\frac{m(m-1)}{2}$ linear independent left eigenvectors of $J_{g,A}J_{g,A}^T$ in the form $w=\textnormal{vec}^T(R_0X^T)$, with corresponding zero eigenvalue. 

Finally, we make the last observation. Since $X$ is a fixed point, i.e., \eqref{eq:vectorized_for_jacobian_Ym} is satisfied for $x=\text{vec}(X^T)$, then for any matrix $R_0 \in \mathbb{R}^{m \times m}$, it holds that $R_0X^T(A-D)^T=0$. Vectorizing this expression, we get 
\begin{align}
    & \text{vec}(R_0X^T(A-D)^T)=\\ \nonumber
    &((A-D)\otimes I_m) \text{vec}(R_0X^T)=0.
\end{align} Thus, $\text{vec}^T(R_0X^T)((A-D)^T\otimes I_m) =\text{vec}^T(R_0X^T)((A-D)\otimes I_m)=0$, and consequently, $w=\textnormal{vec}^T(R_0X^T)$ is a left-eigenvector of $J_{g,x}$ with corresponding zero eigenvalue. Thus, matrices $J_{g,A}$ and $J_{g,x}$ share a null space of dimension $\frac{m(m-1)}{2}$. Thus it follows $\text{rank}(J_g) \leq nm - \frac{m(m-1)}{2}$.

\end{document}